\documentclass[12pt, reqno]{amsart}

\usepackage{latexsym, amsfonts, amsmath, amsthm, amscd} 
\newcommand{\tsig}{\widetilde{\sigma}}

\newcommand{\bil}[2]{\cP(#1,#2)}

\newcommand{\Dodd}{\DD\odd} 
\newcommand{\Deven}{\DD\even}

\newcommand{\even}{^{even}} 
\newcommand{\odd}{^{odd}}

\newcommand{\tsum}{\textstyle\sum}

\newcommand{\boxt}{\Tox\kern -6.3pt\raise .55pt
    \hbox{$\scriptstyle{\times}$}}

\newcommand{\bq}{{\mathbf q}}

\newcommand{\bp}{{\mathbf p}}

\newcommand{\cdott}{\cdot,\cdot}

\newcommand{\D}{{\mathfrak D}}

\newcommand{\RO}[1][]{R^{#1}(\OO)} 
\newcommand{\RM}[1][]{R^{#1}(\M)} 
\newcommand{\RtO}[1][]{R^{#1}(\tO)}

\newcommand{\tO}{\widetilde{\OO}}

\newcommand{\halfN}{\half{\mathbb N}}

\newcommand{\ltimes}{\vbox to 5.4pt{\leaders\vrule\vfil}\kern
  -5pt\times}

\newcommand{\llongrightarrow}{\relbar\joinrel\longrightarrow}

\newcommand{\mmaprightd}[1]{\smash{\mathop
   {\llongrightarrow}\limits^{#1}}}

\newtheorem{thm}{Theorem}[section]
\newtheorem{lem}[thm]{Lemma}

\newtheorem{prop}[thm]{Proposition}
\newtheorem{cor}[thm]{Corollary}

\theoremstyle{definition}
\newtheorem{definition}[thm]{Definition}

\theoremstyle{remark}
\newtheorem{remark}[thm]{Remark}

\numberwithin{equation}{section}

\newcommand{\mapdown}[1]{\Big\downarrow\rlap{$\vcenter
{\hbox{$\scriptstyle#1$}}$}}

\newcommand{\inv}{^{-1}}

\newcommand{\half}[1][1]{\frac{#1}{2}}
\newcommand{\thalf}{\textstyle\frac{1}{2}}

          \newcommand{\g}{{\mathfrak g}}

\newcommand{\C}{{\mathbb C}}
\newcommand{\Z}{{\mathbb Z}}

\newcommand{\bbN}{{\mathbb N}} \newcommand{\bbT}{{\mathcal T}}

\newcommand{\al}{\alpha}
\newcommand{\be}{\beta}
\newcommand{\ga}{\gamma}

\newcommand{\kap}{\kappa}
 \newcommand{\La}{\Lambda}

\newcommand{\sig}{\sigma}

                                     \newcommand{\A}{{\mathcal A}}
                                      
\newcommand{\cC}{{\mathcal C}}

                                             \newcommand{\I}{{\mathcal I}}

                                                \newcommand{\M}{{\mathcal M}}
 
                                          \newcommand{\OO}{{\mathcal O}} 
\newcommand{\cP}{{\mathcal P}} 
                                            \newcommand{\Q}{{\mathcal Q}} 
\newcommand{\cS}{{\mathcal S}}

\newcommand{\gr}{\mathop{\mathrm{gr}}\nolimits}

\newcommand{\End}{\mathop{\mathrm{End}}\nolimits}

\newcommand{\Ker}{\mathop{\mathrm{Ker}}\nolimits}

\newcommand{\gog}{\g\oplus\g}

\newcommand{\Ug}[1][]{{\mathcal U}_{#1}({\mathfrak g})} 
\newcommand{\Zg}{{\mathcal Z}({\mathfrak g})} 
  
\newcommand{\RR}[1][]{{\mathcal R^{#1}}}

\newcommand{\Killg}[1]{({#1})_{\g}}

\newcommand{\GtS}{G\times\cS}

\newcommand{\DD}[1][]{{\mathcal D^{#1}}} 
\newcommand{\DS}[1][]{{\DD^{\cS}_{#1}}}

\begin{document}   

\title 
{From Dixmier algebras to  Star Products} 
\author{Ranee Brylinski}
\address{Department of Mathematics,
        Penn State University, University Park 16802}
\email{rkb@math.psu.edu}
\urladdr{www.math.psu.edu/rkb}


\begin{abstract}
Let $\M$ be a Galois cover of a nilpotent coadjoint orbit of a complex
semisimple Lie group.  We define the notion of a \emph{perfect}  Dixmier
algebra for $\M$ and show how this produces a graded (non-local)
equivariant star  product on  $\M$ with several very nice properties.
This is part of a larger program we have been developing for working out
the orbit method for nilpotent orbits.
\end{abstract}
 
\maketitle      
 
\section{Introduction}  
\label{sec1} 
This note is a companion piece to \cite{me:2k}, and a bridge between
the   approaches  to quantization of nilpotent orbit covers   in
\cite{me:2k} and its sequel \cite{me:2kpos}.  See the introduction to
\cite{me:2k} for background, references and motivations.

The purpose of  this note is to set out some basic relations between the
representation theoretic notion of Dixmier algebras and the Poisson geometric
notion of star products.

\section{Perfect Dixmier algebras}  
\label{sec2}
Let  $\OO$   be a   complex
nilpotent orbit in  some complex semisimple Lie algebra $\g$. We
assume $\OO$ spans $\g$. Let $G^{sc}$ be a simply-connected complex
Lie group  with Lie algebra $\g$.

Let $\kap:\M\to\OO$ be a  covering where $\M$ is
connected. Then the following  geometric structure lifts from $\OO$ to
$\M$:  the adjoint action of $G^{sc}$,
the KKS symplectic form, the Hamiltonian functions
$\phi^x$, $x\in\g$ defined by $\phi^x(y)=\Killg{x,y}$, and the square of
the Euler dilation action of $\C^*$.
In fact,  $G$ acts on $\M$ where $G$ is the quotient of $G^{sc}$ by the
(finite central) subgroup which fixes $\M$.

Consider the algebra 
\begin{equation}\label{eq:RR=} 
\RR=\RM
\end{equation}  
of regular functions on $\M$. Then the Euler grading
$\RR=\oplus_{j\in\halfN}\RR[j]$  makes $\RR$ into a  graded  Poisson 
algebra  in the sense of 
\cite[Definition 2.2.1]{me:2k} with $\phi^x\in\RR[1]$,
for $x\in\g$, and $\{\phi^x,\phi^y\}=\phi^{[x,y]}$.
Also $\RR$ is a superalgebra with even and odd parts
$\RR\even=\oplus_{j\in\bbN}\RR[j]$ and
$\RR\odd=\oplus_{j\in\bbN+\half}\RR[j]$.
Let $\al$ be the algebra automorphism 
of  $\RR$ such that $\al(f)=i^{2j}f$ if $f\in\RR[j]$.
Notice $\al^4=1$. 

Suppose $m\in\M$ covers $e\in\OO$. Then $\M$ is a Galois  cover
of $\OO$  if and only if $G^m$ is normal in $G^e$. In this event,
$\cS=G^e/G^m$ is the Galois group , $\cS$ acts on $\M$ by 
symplectic automorphisms, and our grading of $\RR$ is $\cS$-invariant.
The universal cover of $\OO$ is always Galois.

Now we strengthen the usual definition of Dixmier algebra
(Vogan, McGovern) by adding three new axioms in 
(V), (VI) and (VII). We allowed redundancies in
our axioms   in order to make them as explicit as possible. 
\begin{definition}\label{def:dix} 
Assume the cover $\M$ of $\OO$ is  Galois  with Galois group $\cS$.
A \emph{perfect Dixmier algebra} for $\M$ is a noncommutative 
algebra  $\DD$ together with the following data:
\begin{enumerate} 
\item[(I)] An increasing algebra  filtration 
$\DD=\cup_{j\in\halfN}\,\DD_j$ such that
$[\DD_j,\DD_k]\subseteq\DD_{j+k-1}$ for all $j,k\in\halfN$.
\item[(II)] A representation of $\cS$ on $\DD$ by filtered
algebra automorphisms.
\item[(III)] A Lie algebra embedding $\psi:\g\to\DS[1]$ such that the 
representation of $\g$ on $\DD$ given by the derivations
$a\mapsto\psi^xa-a\psi^x$ is locally finite and 
exponentiates to a  representation  of $G$ on  $\DD$ by algebra
automorphisms.
\item[(IV)] An   $\cS$-equivariant  graded Poisson algebra isomorphism
$\ga:\gr\DD\to\RR$ such that $\ga(\bp_1(\psi^x))=\phi^x$
where $\bp_1:\DD_1\to\gr_1\DD$ is the natural projection.
\item[(V)] An   $\cS$-invariant filtered algebra anti-automorphism 
$\beta$ of  $\DD$   such that 
(a)  $\beta(\psi^x)=\psi^{-x}$, (b) $\be$ induces
$\ga\inv\al\ga$ on $\gr\DD$ and (c) $\be^4=1$.
\end{enumerate}
We impose two further axioms. 
To state these, we notice two useful consequences of  (I)-(V).
First there is a unique $G$-linear map $T:\DD\to\C$ such that $T(1)=1$.
Second, $\DD$ has become a superalgebra with $(\GtS)$-invariant
filtered algebra $\Z_2$-grading
$\DD=\Deven\oplus\Dodd$ where the summands are the 
$\pm 1$-eigenspaces of  $\beta^2$. An element $a\in\DD$ is 
\emph{superhomogeneous} if $a\in\Deven$, in which case $|a|=0$,
or $a\in\Dodd$, in which case $|a|=1$.  Now we require that 
\begin{enumerate}
\item[(VI)]  $T$ is a supertrace.
\item[(VII)] The   $G$-invariant supersymmetric bilinear pairing 
$\bil{a}{b}=T(ab)$ is non-degenerate on $\DD_j$ for each $j\in\halfN$.
\end{enumerate}
\end{definition}
In (VI), $T$ is a supertrace means that if $a$ and $b$ are 
superhomogeneous  then 
\begin{equation}\label{eq:Tab} 
T(ab)=(-1)^{|a||b|}T(ba)  
\end{equation}
when  $a$ and $b$ have same parity, while
$T(ab)=0$ when  $a$ and $b$ have different parity. 
Axioms (IV) and (V)  guarantee  that  $T$ vanishes on $\Dodd$ and so
axiom (VI) amounts to (\ref{eq:Tab}).
 
To make sense of  (VII), we notice if $L\subseteq\DD$ is any
$\beta^2$-stable subspace then we have a  notion of the
$\cP$-orthogonal subspace $L^\perp$ (since the right and left orthogonal
subspaces coincide). We then say $\cP$ is \emph{non-degenerate} on   
$L$ if  $L^\perp\cap L=0$.
The pairing $\cP$ is \emph{supersymmetric} in the sense that
$\Deven$ and $\Dodd$ are $\cP$-orthogonal,
$\cP$ is symmetric on $\Deven$ and  $\cP$ is anti-symmetric on $\Dodd$.
Notice that (VII) is much stronger than  saying $\cP$ is non-degenerate on
$\DD$.

We often speak of $\DD$ as the perfect Dixmier algebra, with the
additional  data being understood.   See \cite[\S7-8]{me:2k} for examples.

Here are some consequences of the axioms.
First  $\RR$ is $\bbN$-graded if and only if $\DD=\Deven$.
Indeed, $\RR$ is $\bbN$-graded $\Leftrightarrow$
$\al^2=1$ $\Leftrightarrow$ $\beta^2=1$, where the last
equivalence follows by axiom (V)(b).
An instance where $\RR$ is $\bbN$-graded is when $\M=\OO$.
Second, $\DS$ is a perfect Dixmier algebra for $\OO$.  
This follows as $\RR^\cS=\RO$ and
all the  Dixmier algebra  data for $\M$ is $\cS$-equivariant. 

Third, axiom (IV) provides a ``symbol calculus" for
$\DD$ with values in $\RR$. 
If $a\in\DD_j$ then the \emph{$\ga$-symbol of a of order $j$}
is the image of $a$ under the map
\[\DD_j\to\DD_j/\DD_{j-\half}\to\RR[j]\]

Fourth, let us extend $\psi:\g\to\DS[1]$ to an algebra homomorphism
\begin{equation}\label{eq:psi:Ug->DS} 
\psi:\Ug\to\DS  
\end{equation}
and let $J$ be the kernel of (\ref{eq:psi:Ug->DS}). Then
$J\cap\Zg$ is a maximal ideal of $\Zg$ where $\Zg$ is the center
of $\Ug$. This follows since by axioms (III) and (IV) the vector
spaces $\D$, $\gr\D$ and $\RtO$ are all $\g$-isomorphic and 
so in particular $\D^G=\C$.

Now (V)(a) says the following square is commutative: 
\begin{equation}\label{eq:tau_beta_sq} 
\begin{array}{ccc}
\Ug&\mmaprightd{\psi}&\DD\\[8pt]
\mapdown{\tau}&&\mapdown{\beta}\\[8pt]
\Ug&\mmaprightd{\psi}&\DD 
\end{array}   
\end{equation}
where $\tau$ is the principal anti-automorphism of $\Ug$.
Consequently  $J$  is a $\tau$-stable $2$-sided ideal in $\Ug$. 
Furthermore $J$ is a completely prime primitive ideal in $\Ug$.
($\Ug/J$ is a subalgebra of $\DD$ and so has no zero-divisors.
This means $J$ is completely prime. But also $J\neq\Ug$ and $J$
contains a maximal ideal of the center of $\Ug$ and so, by a result
of Dixmier, $J$ is primitive.)
 
\section{Simplicity of   $\DD$}  
\label{sec3}

\begin{lem}\label{lem:simple} 
Suppose axioms \textup{(I)-(VI)} are satisfied. Let $\cC$ be some
$\beta^2$-stable subalgebra of $\DD$.  
Then the following two conditions are equivalent:
\begin{itemize}
\item[\rm(i)] $\cC$ is a simple ring
\item[\rm(ii)]  $\cP$ is non-degenerate on $\cC$
\end{itemize}
\end{lem}
\begin{proof} 
(i)$\Rightarrow$(ii):
$\cC$ is  simple implies    $\cC^\perp\cap\cC=0$ since
$\cC^\perp\cap\cC$ is a $2$-sided ideal in $\cC$ 
which does not contain $1$.
\newline
(ii)$\Rightarrow$(i):
Let $\I$ be a non-zero two-sided ideal in $\DD'$.
Pick $a\in\I$  with $a\neq 0$. Then (ii) implies  there exists 
$b\in\cC$ such that $T(ab)=1$. It follows that $ab=c+1$ 
where $c$ lies in  $\Ker T\cap\cC$. 
So $\I$ contains $c+1$. But then  $\I$ contains the 
$G$-subrepresentation generated by $c+1$.
Since $\Ker\,T$ contains no non-zero $G$-invariants it follows 
(by completely  reducibility of $\cC$ as a  $G$-representation), 
that $\I$ contains  both $c$ and $1$. Thus   $\cC$ is simple.
\end{proof}

\begin{prop}\label{prop:simple_iff} 
Suppose we are in the situation of Definition \textup{\ref{def:dix}}
and axioms \textup{(I)-(VI)} are satisfied. Then
$\DD$ is a simple ring if and only if $\DS$ is a simple ring.
\end{prop}
\begin{proof} 
Suppose $\DD$ is simple. Then $\cP$ is non-degenerate on $\DD$
(by Lemma \ref{lem:simple}) and hence 
(since  $\cP$ is $\cS$-invariant) $\cP$ is non-degenerate on $\DS$.
Then  (by Lemma \ref{lem:simple} again) $\DS$ is simple.

Conversely, assume $\DS$ is simple. 
Let $\I$ be a non-zero $2$-sided ideal in $\DD$.
To show $\DD$ is simple, 
it suffices to show that $\I^{\cS}=\I\cap\DS$ is non-zero.
Let $a\in\I$. Consider 
$b=(n!)\inv\sum_{\sig\in S_n} (s_{\sig_1}a)\cdots(s_{\sig_n}a)$
where $s_1,\dots,s_n$ is some listing of the elements of the 
Galois group  $\cS$.
Clearly $b\in\I^{\cS}$. But also we can see using axiom (IV) that 
$b\neq 0$.
Indeed, let $j$ be the filtration order of $a$ in $\DD$ and let 
$\phi\in\RR[j]$ be the  $\ga$-symbol  of order $j$ of $a$. Then
$b$ lies in $\DD_{jn}$ and the  $\ga$-symbol  of order $jn$ of $b$ is
$(s_{\sig_1}\phi)\cdots(s_{\sig_n}\phi)$. 
This product is non-zero and so $b$  must be non-zero.
\end{proof}

\begin{cor}\label{cor:simple} 
If $\DD$ is a perfect Dixmier algebra then
$\DD$ and $\DS$ are both simple rings. 
\end{cor}

\begin{cor}\label{cor:Ug_DS} 
If $\DD$ is a perfect Dixmier algebra and \textup{(\ref{eq:psi:Ug->DS})} 
is surjective,  then the kernel $J$ is a maximal ideal in $\Ug$.
\end{cor}

\begin{remark}\label{rem:Ug/J=simple}  
Ideally the axioms for a perfect Dixmier algebra  should
automatically imply that the image of (\ref{eq:psi:Ug->DS}) is simple,
i.e., that $J$ is maximal. 
Our axioms (I)-(VII) may not do this, but in any event our enriched axiom
set  (see Remark \ref{rem:pos})   accomplishes this.
\end{remark}

\section{The noncommutative   $\circ$ product}  
\label{sec4}

By axiom (VII) in Definition \ref{def:dix}, 
we have a unique   $\cP$-orthogonal  decomposition
\begin{equation}\label{eq:DD=} 
\DD=\bigoplus_{j\in\halfN}\DD^j  
\end{equation}
such that
$\DD_k=\bigoplus_{j=0}^k\DD^j$.
Each  space $\DD^j$ is $(\GtS)$-stable. It is easy to see that
\[\cP(a,b)=\cP(\beta(b),\beta(a))\] 
for all $a,b\in\DD$. Consequently 
$\DD^j$ is $\beta$-stable and using axiom (V)(b) we find that
$\beta$ acts on $\DD^j$ by multiplication by $i^{2j}$.
Then   $\Deven=\bigoplus_{j\in\bbN}\DD^j$,
$\Dodd=\bigoplus_{j\in\bbN+\half}\DD^j$ and $T$ is the
orthogonal projection of $\DD$ onto $\DD^0=\C$.
Clearly there is  a   unique linear map 
\begin{equation}\label{eq:bq} 
\bq:\RR\to\DD
\end{equation}  
such that $\bq$ lifts $\ga\inv:\RR\to\gr\DD$ and 
$\bq(\RR[j])=\DD^j$ for  all   $j\in\halfN$. 
Then $\bq$ is a  $(\GtS)$-linear vector space isomorphism,
and we have $\psi^x=\bq\inv\phi^x$ and  $\beta=\bq\al\bq\inv$.
Now we define, for all $\phi,\psi\in\RR$,
\begin{equation}\label{eq:circ=} 
\phi\circ\psi=\bq\inv((\bq \phi)(\bq\psi))
\end{equation}
\begin{prop}\label{prop:circ} 
Assume we have a perfect Dixmier algebra $\DD$ for  $\M$ with $\DD^j$ 
and $\bq$ defined as  above. Then $\circ$ is a $(\GtS)$-invariant 
associative product on $\RR$ and so $(\RR,+,\circ)$ is an
associative noncommutative algebra.  Then
\begin{equation}\label{eq:circj-k} 
\RR[j]\circ\RR[k]\subseteq
\RR[j+k]\oplus\RR[j+k-1]\oplus\cdots\oplus\RR[|j-k|]
\end{equation}
where  $j,k\in\halfN$. Suppose $\phi\in\RR[j]$ and $\psi\in\RR[k]$ so that
$\phi\circ\psi=\sum_pC_p(\phi,\psi)$ where
$C_p(\phi,\psi)$ lies in $\RR[j+k-p]$. Then
\begin{eqnarray}
\phi\circ\psi&\equiv&\phi\psi+\half\{\phi,\psi\}\mod\RR[\le j+k-2]
\label{eq:circ_01}\\[3pt]
C_p(\phi,\psi)&=&(-1)^p\, C_p(\psi,\phi)
\label{eq:circ_par}
\end{eqnarray}
\end{prop}
\begin{proof} 
The second sentence is clear. Now (\ref{eq:circj-k}) is equivalent to
\begin{equation}\label{eq:DDj-k} 
\DD[j]\DD[k]\subseteq
\DD[j+k]\oplus\DD[j+k-1]\oplus\cdots\oplus\DD[|j-k|]
\end{equation}
We have 
$\DD^j\DD^k\subseteq\bigoplus_{i\in\bbN}\DD^{j+k-i}$ 
since $\DD$ is a superalgebra.
Because of axiom (VII), showing  (\ref{eq:DDj-k}) reduces to
showing that $\DD^j\DD^k$ is orthogonal to $\DD^s$ if $s<|j-k|$.
So suppose $\DD^j\DD^k$ is not orthogonal to $\DD^s$. Then
there exist $a\in\DD^j$, $b\in\DD^k$ and $c\in\DD^s$ such that
$T(abc)=1$. Then $bc$ has a component in $\DD^j$ and so
$k+s\ge j$. But also $T(bca)=\pm1$ (since $T$ is a supertrace) and so
similarly $s+j\ge k$. Hence $s\ge|j-k|$. Now for 
$\phi\in\RR[j]$ and $\psi\in\RR[k]$   we can write
\begin{equation}\label{eq:Cp} 
\phi\circ\psi=\sum_{p=0}^{2\min{(j,k)}}C_p(\phi,\psi)
\end{equation}
where $C_p(\phi,\psi)\in\RR[j+k-p]$.

Now axiom  (IV) implies that 
$C_0(\phi,\psi)=\phi\psi$ and  
$C_1(\phi,\psi)-C_1(\psi,\phi)=\{\phi,\psi\}$.
But also  axiom (V) implies that the map $\al:\RR\to\RR$ (which is an 
algebra automorphism with respect to the ordinary product)
is an algebra anti-automorphism with respect to $\circ$.  Thus 
\[\al(\phi\circ\psi)=(\al\psi)\circ(\al\phi)=i^{2j+2k}\,\psi\circ\phi\]
Then $i^{-2p}C_p(\phi,\psi)=C_p(\psi,\phi)$.
This proves (\ref{eq:circ_par}). Then in particular
$C_1(\phi,\psi)=-C_1(\psi,\phi)$. So $C_1(\phi,\psi)=\half\{\phi,\psi\}$
and  we get (\ref{eq:circ_01}).
\end{proof}

We can think of $\bq:\RR\to\DD$ a  ``quantization map" as in
\cite[\S8.1-8.2]{me:2k}. In particular we have 
\begin{cor}\label{cor:exact} 
If $\phi\in\RR[1]$, for instance if $\phi=\phi^x$ where $x\in\g$, then
for all $\psi\in\RR$  we have 
$\{\phi,\psi\}=\phi\circ\psi-\psi\circ\phi$.
\end{cor}
\begin{proof} 
Identical to the proof of \cite[Corollary  8.2.3]{me:2k}.
\end{proof}

Proposition \ref{prop:circ} implies in particular that
$\RR$, equipped with its $\circ$ product, 
\emph{becomes the perfect Dixmier algebra}!
The data on $\DD$ required by the axioms corresponds under $\bq$ to
data that  exist  from the beginning on $\RR$:

(I)  the grading on $\RR$ gives rise to the filtration by subspaces
$\RR[\le j]$,

(II) the $\psi^x$ correspond to the momentum functions $\phi^x$ of 
the Hamiltonian $\g$-symmetry,

(III)  the Galois group $\cS$ already acts on $\RR$,

(IV) $\ga$ corresponds to the identity map,

(V)  $\beta$ corresponds to $\al$.

\noindent Notice that the even and odd parts of $\DD$ correspond to the 
even and odd  parts of  $\RR$. 
Also  $T$ corresponds to the projection  
\begin{equation}\label{eq:bbT} 
\bbT:\RR\to\RR[0]=\C  
\end{equation}
defined by the Euler grading. 
So $\bbT(\phi)$ is just the \emph{constant term} of $\phi$. 
The nondegenerate bilinear  pairing $\Q$ on $\RR$ corresponding to $\cP$ is
given by $\Q(\phi,\psi)=\bbT(\phi\circ\psi)$. This pairing is 
\emph{graded}  supersymmetric    in the sense that
$\Q$ pairs $\RR[j]$ with $\RR[k]$ trivially if $j\neq k$.

In this   approach, the ``new" axioms (VI) and (VII) have played a crucial 
role. Also we  have gained a lot more    structure on the Dixmier algebra.
In particular  the $\circ$ product ``breaks off" after
degree $|j-k|$ in (\ref{eq:circj-k}); this is a Clebsh-Gordan type 
phenomenon.
 
Thus the problem of finding  a Dixmier algebra for $\M$ can be
reformulated as  the problem of constructing a suitable product $\circ$ 
on  $\M$. To formalize this we make 
\begin{definition}\label{def:dixp} 
Assume $\M$ is a Galois cover of $\OO$ with Galois group $\cS$.
A \emph{perfect Dixmier product}  on $\RR$ is a $(\GtS)$-invariant
associative  noncommutative product  $\circ$ 
satisfying (\ref{eq:circj-k}), (\ref{eq:circ_01}) and (\ref{eq:circ_par}) 
such  that the bilinear pairing 
$\Q(\phi,\psi)=\bbT(\phi\circ\psi)$ on $\RR$ is
graded supersymmetric and non-degenerate.
\end{definition}

A  perfect  Dixmier product  makes $\RR$ into  a  filtered superalgebra
which  then, together with  the Hamiltonian functions $\phi^x$, $x\in\g$,
and   $\al$,  is a perfect Dixmier algebra  for $\M$ in the sense  of
Definition \ref{def:dix}. Conversely, we have shown that a  perfect
Dixmier algebra for $\M$ yields a  perfect Dixmier product on $\RR$.

\begin{remark}\label{rem:pos}
We can enrich our axiom set to get the notions for $\M$ of a 
\emph{positive Dixmier algebra} and a \emph{positive Dixmier product}.
The extra axioms require   lifting  the complex  conjugation
automorphism $\sig$ of $\OO$  (induced by the Cartan involution  of $\g$)
to an antiholomorphic automorphism 
$\tsig$ of $\M$ (of order $2$ or $4$) such that   
(i) $\tsig$ induces a $\C$-antilinear $\circ$-algebra automorphism
of $\RR$ and (ii) the Hermitian pairing 
$(\phi|\psi)=\bbT(\phi\circ\psi^{\tsig})$ is  positive-definite.
($\bbT$ being a supertrace is equivalent to  this pairing being Hermitian.)
We   develop  this in  \cite{me:2kpos} in the context of star products.
\end{remark}

\section{Graded star products}  
\label{sec5}

Let $\A=\oplus_{j\in\halfN}\A^j$  
be a graded Poisson algebra as in \cite[Definition 2.2.1]{me:2k}.
Assume $\A^0=\C$.
Then a \emph{graded star product} (with parity) on     
$\A$ is a product
$\star$ on $\A[t]$ which makes $\A[t]$  into an   associative algebra
over  $\C[t]$ such that, for $\phi,\psi\in\A$, the series 
$\phi\star\psi=\tsum_{p=0}^\infty C_p(\phi,\psi)t^p$ satisfies:

(i) $C_0(\phi,\psi)=\phi\psi$ 

(ii) $C_1(\phi,\psi)=\thalf\{\phi,\psi\}t$ 

(iii)  $C_p(\phi,\psi)=(-1)^p C_p(\psi,\phi)$ 

(iv) $C_p(\phi,\psi)\in\A^{j+k-p}$ when $\phi\in\A^j$ and 
$\psi\in\A^k$

\noindent We do \emph{not} require that $\star$ is bidifferential, i.e. that the
operators $C_p(\cdott)$ are bidifferential.

A graded star product $\star$ on $\A$ is   specializes at $t=1$ to give
a noncommutative product $\circ$ on $\A$. 
Clearly $\circ$ uniquely determines $\star$. 
Let $T:\A\to\A^0=\C$ be the projection defined by the grading.
Then we have a bilinear pairing 
$\Q:\A\times\A\to\C$ defined by 
\begin{equation}\label{eq:Q} 
\Q(\phi,\psi)=T(\phi\circ\psi)
\end{equation}
Notice that  $\Q(\phi,\psi)=C_{2k}(\phi,\psi)$ if $\phi,\psi\in\A^k$
and so  the parity axiom (iii) implies that $\Q$ is symmetric if 
$k\in\bbN$  or anti-symmetric if $k\in\bbN+\half$.

\begin{definition}\label{def:perfect} 
We say $\star$  is  \emph{orthogonally graded} if
$\A^j$ and $\A^k$ are $\Q$-orthogonal when $j\neq k$.
We say $\star$  is  \emph{perfectly graded} if also
$\Q$ is non-degenerate on $\A^j$ for each $j$.
\end{definition}

Comparing with the proof of Proposition \ref{prop:circ}, we find
\begin{lem}\label{lem:stuff} 
If $\star$  is orthogonally graded  then 
$T$ is a supertrace on $\A$ with respect to $\circ$. If
$\star$  is perfectly graded then 
\begin{equation}\label{eq:Aj*Ak} 
\A^j\star\A^k\subset
\A^{j+k}\oplus t\A^{j+k-1}\oplus\cdots\oplus t^{2\,\min(j,k)}\A^{|j-k|} 
\end{equation}
\end{lem}

Suppose  we have Hamiltonian symmetry $\phi:\g\to\A^1$,
$x\mapsto\phi^x$, as in \cite[Definition 3.1.1]{me:2k}
Put $[\phi,\psi]_\star=\phi\star\psi-\psi\star\phi$.
We say $\star$ is
\emph{$\g$-covariant} if 
$[\phi^x,\phi^y]_\star=t\phi^{[x,y]}$ for all $x,y\in\g$.
We say $\star$ is \emph{exactly $\g$-invariant}
(or \emph{strongly $\g$-invariant}) if  we have the 
stronger property:
\begin{equation}\label{eq:exact} 
[\phi^x,\psi]_\star=t\{\phi^x,\psi\}
\end{equation}
for all $x\in\g$ and $\psi\in\A$.
Exact $\g$-invariance implies ordinary $G$-invariance, i.e.,
$(g\cdot\psi_1)\star(g\cdot\psi_2)=g\cdot(\psi_1\star\psi_2)$
where $G$ acts on $\A[t]$ by $g\cdot(\psi t^p)=(g\cdot\psi)t^p$.

\begin{lem}\label{lem:exact} 
If $\star$ is an orthogonally  graded star product on $\A$ then
$\star$ is   exactly $\A^1$-invariant.
\end{lem}
\begin{proof} 
If $\phi\in\A^1$ and $\psi\in\A$ then 
$\phi\star\psi=\phi\psi+\half\{\phi,\psi\}t+t^2C_2(\phi,\psi)$.
Then $[\phi^x,\psi]_\star=t\{\phi^x,\psi\}$ because of 
the parity  axiom (iii). 
\end{proof}

\begin{lem}\label{lem:circ-star} 
Suppose $\circ$ is a perfect Dixmier product on $\RR$.
Then $\circ$ is the specialization  at $t=1$
of a unique graded star product $\star$ on $\RR$. 
Moreover $\star$ is perfectly graded and $\cS$-invariant.

Conversely, suppose $\star$ is a  perfectly graded, $\cS$-invariant 
star product on $\RR$. 
Then  the  specialization at $t=1$ of $\star$ is a 
perfect Dixmier product on $\RR$.
\end{lem}
\begin{proof} 
Given $\circ$, we define  a product $\star$ on $\RR{}[t]$ as follows:
$\star$ is  $\C[t]$-bilinear and if
$\phi\in\RR[j]$ and $\psi\in\RR[k]$ with
$\phi\circ\psi=\sum_{i=|j-k|}^{j+k}\pi_{i}$ where $\pi_i\in\RR[i]$ then
$\phi\star\psi=\sum_{i=|j-k|}^{j+k}\pi_{i}\,t^{j+k-i}$.  
This is the only possible way to   extend $\circ$ to a graded star
product. The properties of $\circ$ imply that
$\star$ is in fact a perfectly graded, $\cS$-invariant star product.
The converse is clear.
\end{proof}

\subsection{The operators \boldmath $\La^x$}  
\label{ss7_Lax} 
Here is an important consequence of perfectness.
\begin{prop}\label{prop:Lax} 
Suppose we have a  perfect Dixmier product $\circ$ on $\RR$.
For $x\in\g$ and any $\psi\in\RR$ we have
\begin{equation}\label{eq:phix*} 
\phi^x\circ\psi=\phi^x\psi+\half\{\phi^x,\psi\}+\La^x(\psi)
\end{equation}
where   $\La^x:\RR\to \RR$ are linear operators. These satisfy  
\begin{itemize}
\item[\rm(i)] $\La^x$ is the $\Q$-adjoint of ordinary multiplication by 
$\phi^x$.
\item[\rm(ii)] If $x\neq 0$ and $j$ is positive, then $\La^x$ is   non-zero
somewhere on $\RR[j]$.
\item[\rm(iii)] $\La^x$ is graded of degree $-1$, i.e.,
$\La^x(\RR[j])\subseteq\RR[j-1]$.
\item[\rm(iv)] The operators $\La^x$ commute, i.e., 
$[\La^x,\La^y]=0$ for  all $x,y\in\g$.
\item[\rm(v)] The operators $\La^x$ transform in the adjoint
representation of $\g$, i.e., $[\Phi^x,\La^y]=\La^{[x,y]}$
where $\Phi^x=\{\phi^x,\cdot\}$.
\item[\rm(vi)] The operators $\La^x$ commute with the 
$\cS$-action on  $\RR$.
\end{itemize}
\end{prop}  
\begin{proof} 
Same as the proof of \cite[Corollary 8.4.1]{me:2k}.
\end{proof}

If we identify $\RR$ with $\DD$ via $\bq$, then we get the
representation   
\begin{equation}\label{eq:gog_RtO} 
\Pi:\gog\to\End_{\cS}\RR,\qquad 
\Pi^{(x,y)}(\psi)=\phi^x\circ\psi-\psi\circ\phi^y
\end{equation} 
Then    Proposition \ref{prop:Lax} gives
\begin{cor}\label{cor:phi+La} 
For $x\in\g$ we have $\Pi^{(x,x)}=\eta^x$ and 
$\Pi^{(x,-x)}=2\phi^x+2\La^x$.
\end{cor}

\bibliographystyle{plain}

\end{document}